\documentclass[10pt,a4paper]{amsart}
\usepackage{amsfonts,amsmath,amsthm,amsopn}

\newcommand{\R}{\mathbb{R}}
\renewcommand{\ge}{\varepsilon}
\newcommand{\psa}{{\Psi_{\alpha}}}
\renewcommand{\O}{{\Omega}}

\DeclareMathOperator{\supp}{supp}

\newtheorem{theorem}{Theorem} 
\newtheorem{lemma}[theorem]{Lemma}

\newtheorem{proposition}[theorem]{Proposition}
\newtheorem{corollary}[theorem]{Corollary}

\theoremstyle{definition}
\newtheorem{remark}[theorem]{Remark}

\renewcommand{\a}{\alpha}

\date{\today}
\begin{document}
\title[MULTIPLE SOLUTIONS FOR A HENON-LIKE EQUATION ON THE ANNULUS]{MULTIPLE SOLUTIONS FOR A HENON-LIKE EQUATION ON THE ANNULUS}
\author{Marta Calanchi}
\address{Dipartimento di Matematica,
  Universit\`a degli Studi di Milano, via C. Saldini 50, 20133 Milano,
  Italy}
\email{Marta.Calanchi@mat.unimi.it}
\author{Simone Secchi}
\address{Dipartimento di Matematica e Applicazioni, Universit\`a di Milano--Bicocca,
via R. Cozzi 53, I-20125 Milano} \email{Simone.Secchi@unimib.it}
\author{Elide Terraneo}
\address{Dipartimento di Matematica,
  Universit\`a degli Studi di Milano, via C. Saldini 50, 20133 Milano,
  Italy}
\email{Elide.Terraneo@mat.unimi.it}
\thanks{The second author is partially
  supported by MIUR, national project \textit{Variational Methods and
    Nonlinear Differential Equations}} 

 \begin{abstract}
   For the equation \(
  -\Delta u = \left| |x|-2 \right|^\alpha u^{p-1} \),
  \(1 < |x| < 3\), we prove the existence of two solutions for
  \(\alpha\) large, and of two additional solutions when \(p\) is
  close to the critical Sobolev exponent \(2^*=2N/(N-2)\). A symmetry--breaking phenomenon appears,
  showing that the least--energy solutions cannot be
  radial functions.
\end{abstract}

\maketitle

\section{Introduction}

In this paper we will consider the following problem:
\begin{equation} \label{eq1}
  \begin{cases}
    -\Delta u=\Psi_{\alpha}u^{p-1} &\text{in $\Omega$},\\ u>0
    &\text{in $\Omega$},\\ u=0 &\text{on $\partial \Omega$},
  \end{cases}
\end{equation}
where $\Omega=\{x\in {\mathbb R^N}|\;1<|x|<3\}$ is an annulus in
$\mathbb R^N$, $N\geq 3$, $\alpha>0$, $p>2$ and $\Psi_{\alpha}$ is
the radial function
\[
\Psi_{\alpha}(x)=\left| |x|-2 \right|^{\alpha}.
\]
This equation can be seen as a natural extension to the annular
domain $\Omega$ of the celebrated H\'{e}non equation with
Dirichlet boundary conditions (see \cite{henon,ni})
\begin{equation} \label{eq:henon}
\begin{cases}
-\Delta u = |x|^\alpha |u|^{p-1} &\hbox{for $|x| < 1$} \\
u = 0 &\hbox{if $|x|=1$}.
\end{cases}
\end{equation}
Actually, the weight function $\Psi_\alpha$ reproduces on $\Omega$
a similar qualitative behavior of $|\cdot|^\alpha$ on the unit
ball $B$ of $\R^N$.

A standard compactness argument shows that the infimum
\begin{equation} \label{eq:infimum}
\inf_{\substack{u \in H_0^1(B) \\ u \neq 0}} \frac{\int_B |\nabla
u|^2\ dx}{\left( \int_B |x|^\alpha |u|^p \ dx\right)^{2/p}}
\end{equation}
is achieved for any~$2<p<2^*$ and any~$\alpha>0$. In 1982, Ni
proved in \cite{ni} that the infimum
\begin{equation} \label{eq:infimum-rad}
\inf_{\substack{u \in H_{0,\mathrm{rad}}^1(B) \\ u \neq 0}}
\frac{\int_B |\nabla u|^2\ dx}{\left( \int_B |x|^\alpha |u|^p\ dx
\right)^{2/p}}
\end{equation}
is achieved for any $p \in (2,2^* + \frac{2\alpha}{N-2})$ by a
function in $H_{0,\mathrm{rad}}^1(B)$, the space of radial
$H_0^1(B)$ functions. Thus, radial solutions of \eqref{eq:henon}
exist also for (Sobolev) supercritical exponents $p$. Actually,
radial $H_0^1$ elements show a power--like decay away from the
origin (as a consequence of the Strauss Lemma, see
\cite{AM,strauss}) that combines with the weight $|x|^\alpha$ and
provides the compactness of the embedding $H_{0,\mathrm{rad}}^1(B)
\subset L^p(B)$ for any $2<p<2^*+\frac{2\alpha}{N-2}$.

A natural question is whether any minimizer of \eqref{eq:infimum}
must be radially symmetric in the range $2<p<2^*$ and $\alpha>0$. Since
the weight $| \cdot|^\alpha$ is an increasing function,
neither rearrangement arguments nor the moving plane techniques of
\cite{gnn} can be applied.

Reverting to the case $\alpha >0$, Smets \textit{et al.} proved in
\cite{ssw}
some symmetry--breaking results for \eqref{eq:henon}. They proved
that minimizers of \eqref{eq:infimum} (the so-called
\emph{ground-state
  solutions}, or least energy solutions) cannot be radial, at least
for $\alpha$ large enough. As a consequence, \eqref{eq:henon} has
at least two solutions when $\alpha$ is large (see also
\cite{sw}).

Later on, Serra proved in \cite{serra03} the existence of at least
one non--radial solution to \eqref{eq:henon} in the critical case
$p=2^*$, and in \cite{badialeserra} the authors proved the
existence of more than one solution to the same equation also for
some supercritical values of $p$. These solutions are non-radial
and they are obtained by minimization under suitable symmetry
constraints.

Quite recently, Cao and Peng proved in \cite{caopeng} that, for
$p$ sufficiently close to $2^*$, the ground-state solutions of
\eqref{eq:henon} possess a unique maximum point whose distance
from $\partial B$ tends to zero as $p \to 2^*$.

This kind of result was improved in \cite{peng}, where multibump
solutions for the H\'{e}non equation with almost critical Sobolev
exponent $p$ are found, by means of a finite--dimensional
reduction. These solution are not radial, though they are
invariant under the action of suitable subgroups of $O(N)$, and
they concentrate at boundary points of the unit ball of $\R^N$ as
$p \to 2^*$. The r\^{o}le of $\alpha$ is however a static one.
(For more results for $p\approx 2^*$ see also \cite{PiSe05}).

\medskip

In this paper we will prove that similar phenomena take place for
problem (\ref{eq1}) on the annulus $\Omega$. In
Section~\ref{paralpha}, we present some estimates for the least
energy radial solutions of \eqref{eq1} when $p<2^*$ is kept fixed
but $\alpha \to +\infty$. These will lead us to a first
symmetry--breaking result, stating that for $\alpha$ sufficiently
large there exist at least two solutions of \eqref{eq1}: a global
minimizer of the associated Rayleigh quotient, and a global
minimizer \emph{among radial functions}.

In Section 3, another symmetry--breaking is proved, with $\alpha$
fixed and $p \to 2^*$. To show this phenomenon, we will use a
decomposition lemma in the spirit of P.L.~Lions' concentration and
compactness theory (\cite{PLL}), and inspired by \cite{caopeng}.
It will turn out that global minimizers of the same Rayleigh
quotient concentrate as $p \to 2^*$ at precisely one point of the
boundary $\partial \Omega$, which has two connected components. A
second nonradial solution can then be found in a tricky but
natural way, by minimization over functions that are ``heavier''
on the opposite connected component of $\partial\Omega$.

In Section 4, a third nonradial solution is singled out, by means
of a linking argument. Roughly speaking, the previous nonradial
solutions can be used to build a mountain pass level. In
particular,
 this third solution will not be a local minimizer of
the Rayleigh quotient.

 Section 5 describes the behavior of
ground-state solutions of \eqref{eq1} as $\alpha \to +\infty$ and
$p<2^*$ is kept fixed. Although the conclusion is not as precise
as in the case $p \to 2^*$, we can nevertheless show that a sort
of concentration near the boundary $\partial \Omega$ still
appears.

 We would like to stress that the existence
of non-radial solutions in the annulus in the almost critical case
$p\approx 2^*$ is not by now a surprise. When the weight
disappears, i.e. $\alpha=0$, Brezis and Nirenberg proved in
\cite{BN} that the ground state solution of $-\Delta u = u^p$ in
$H_0^1$ is not a radial function. Indeed, the authors proved that
both a radial and a non-radial (positive) solution arise as $p
\approx 2^*$. Their simple continuation argument can be adapted to
cover our weighted equation. Subgroups of $O(N)$ are used in
\cite{byeon01} for the equation $-\Delta u + u = f(u)$, and some
refined properties of symmetric solutions are proved. We refer to
the bibliography of that paper for more references.

 For more results about asymptotic estimates for
solutions of the H\'{e}non equation with $\alpha$ large, see
\cite{byeonwangI,byeonwang}.

\section{Symmetry breaking for $\alpha$ large }\label{paralpha}

Let $H_{0,\mathrm{rad}}^1(\Omega)$ be the space of radially
symmetric functions of $H_0^1(\Omega)$. With a slight but common
abuse of notation, we will systematically write $u(x)=u(|x|)$ for
$u \in H_{0,\mathrm{rad}}^1(\Omega)$.

Consider the minimization problem
\begin{equation} \label{eq:6} S_{\alpha,p}^{\rm rad} = \inf_{u \in
    H_{0,{\rm rad}}^1 (\Omega) \setminus \{0\}} R_{\alpha,p} (u),
\end{equation}
where
\begin{equation}\label{eq2}
  R_{\alpha,p}(u)= \frac{\int_{\Omega}|\nabla
    u|^2 \, dx}{\left(\int_{\Omega}\psa |u|^p \, dx\right)^\frac
    2p},\qquad u\in H^1_0(\Omega)\setminus \{0\}
\end{equation}
is the \emph{Rayleigh quotient} associated to \eqref{eq1}. It is known
that any minimizers of \eqref{eq:6} can be scaled so as to become
solutions of \eqref{eq1}. Therefore, we will use freely this fact in
the sequel.
\begin{remark}
  Unlike the result of~\cite{ni}, the fact that the annulus $\Omega$
  does not contain the origin implies the existence of a radial
  solution of \eqref{eq1} for any $p>2$. Indeed, the embedding
  $H_{0,\mathrm{rad}}^1(\Omega)\subset L^q(\Omega)$ is compact for
  every $q \geq 1$, and therefore the infimum \eqref{eq:6}
  is achieved by a (radial) function.
\end{remark}

In the next Proposition, we provide an estimate of the
\textit{energy} $S^{\mathrm{rad}}_{\alpha,p} $ as $\alpha \to
\infty$.

\begin{proposition} \label{prop:2}
Let $p>2$. As $\alpha \to \infty$, there exist two constants $C_1$
and $C_2$ depending on $p$ such that
\begin{equation} \label{eq:11}
0< C_1 \leq \frac{S_{\alpha,p}^{\rm rad}}{ \alpha^{1+2/p}} \leq
C_2 < + \infty.
\end{equation}
Moreover, for any $M>2$ it is possible to choose the constants
$C_1$ and $C_2$ independent of $p\in (2,M]$.
\end{proposition}

\begin{proof}
The upper bound $C_2$ can be obtain exactly as in \cite{ssw}: we
fix a positive, radial function $\psi \in C_0^\infty(\Omega)$, and
set $\psi_\alpha (x) =\psi_\alpha (|x|) = \psi (\alpha
(|x|-3+3/\alpha))$. Then
\begin{eqnarray*}
\int_\Omega |\nabla \psi_\alpha|^2 \, dx &=& \omega_{N-1}
\int_{3-\frac{2}{\alpha}}^3
                          \left( \psi'_\alpha(r) \right)^2 r^{N-1} \, dr\\
&=&  \omega_{N-1} \int_{1}^3
                          \alpha^2\psi'(s)^2 \left( \frac{s}{\alpha}
                          +3-\frac{3}{\alpha} \right)^{N-1}
                          \alpha^{-1} \, ds \\
&=& \alpha\ \omega_{N-1} \int_{1}^3 \psi'(s)^2 \left(
\frac{s+3\alpha-3}{\alpha s}
                          \right)^{N-1} \, s^{N-1}\, ds \\
&\leq& 3^{N-1} \alpha \int_{\Omega} | \nabla \psi |^2 \, dx, \quad
\hbox{(since $1 \leq \frac{s+3\alpha-3}{\alpha s}\leq 3$)}
\end{eqnarray*}
and
\begin{equation*}
\int_\Omega \Psi_\alpha  \psi_\alpha^p \, dx \geq \left(1 -
\frac{2}{\alpha} \right)^{\alpha} \alpha^{-1} \int_\Omega \psi^p
\, dx.
\end{equation*}
This proves that $S_{\alpha,p}^{\mathrm{rad}} \leq C(\alpha, p)
\alpha^{1+\frac{2}{p}}$, where
\[
C(\alpha, p) =3^{N-1}\frac{\int_\Omega|\nabla \psi|^2 dx}{(1-\frac
2\alpha)^\frac{2\alpha}{p}(\int_\Omega \psi^p(x)\, dx)^\frac 2p}
\leq C_2 \quad\hbox{for any $p>2$ and $\a>1$}.
\]
To find the lower bound $C_1$, we will perform some scaling. Let
us define the functions $\psi_1 \colon [1,2] \to [1,2]$ and
$\psi_2 \colon [2,3] \to [2,3]$ as follows:
\begin{equation}
\psi_1 (r) = 2- (2- r)^\beta, \quad \psi_2 (r) = 2 +  (r-2)^\beta,
\end{equation}
where $\beta\in (0,1)$ will be chosen hereafter. It is clear that
we can obtain a piecewise $C^1$ homeomorphism  $\psi \colon [1,3]
\to [1,3]$ by gluing $\psi_1$ and $\psi_2$. Now, for any radial
function $u \in H_0^1(\Omega)$, we set $v(\rho) = u(\psi(\rho))$
and compute:
\begin{eqnarray}\label{eq:8}
\int_\Omega |\nabla u|^2 \, dx
&=&\omega_{N-1} \int_{1}^{3} |u'(r)|^2 r^{N-1} \, dr \notag \\
&\geq&\omega_{N-1} \int_{1}^{3} |u'(r)|^2  \, dr \notag \\
&=& \omega_{N-1} \left(\int_{1}^{2} |v'(\rho)|^2
\frac{1}{\psi_1'(\rho)} \, d\rho +
 \int_{2}^3 |v'(\rho)|^2 \frac{1}{\psi_2'(\rho)} \, d\rho\right) \notag \\
&=& \omega_{N-1} \, \frac{1}{\beta}\left( \int_{1}^{2}
|v'(\rho)|^2 (2-\rho)^{1-\beta} \, d\rho \right. \notag \\
&\qquad& \left. {} +\int_{2}^{3}
|v'(\rho)|^2 (\rho-2)^{1-\beta} \, d\rho \right)\notag \\
&=& \omega_{N-1} \, \frac{1}{\beta} \int_{1}^{3}
|v'(\rho)|^2 |\rho-2|^{1-\beta} \, d\rho \notag\\
&\geq&  \omega_{N-1} \, \frac{1}{\beta} \int_{1}^{3} |v'(\rho)|^2
|\rho-2| \, d\rho.
\end{eqnarray}
Choosing $\beta = 1/(\alpha+1)$,
\begin{eqnarray}
\int_\Omega \Psi_\alpha (x) |u(|x|)|^p \, dx
&=& \omega_{N-1} \int_{1}^{3}\Psi_\alpha (r)|u(r)|^p r^{N-1}\; dr\notag \\
&\leq&3^{N-1} \omega_{N-1} \int_{1}^{3}\Psi_\alpha (r)|u(r)|^p\; dr\notag \\
&=&3^{N-1}\omega_{N-1}\left( \int_{1}^{2} \Psi_\alpha
(\psi_1(\rho)) |v(\rho)|^p
\psi_1'(\rho) \, d\rho\right.  \notag \\
&& \qquad + \left. \int_{2}^3 \Psi_\alpha
(\psi_2(\rho))|v(\rho)|^p \psi_2'(\rho) \,
 d\rho\right) \notag \\
&=& 3^{N-1}\omega_{N-1} \beta \int_{1}^3 |v(\rho)|^p \, d \rho.
\label{eq:14}
\end{eqnarray}
Since we are integrating over $\Omega$ and $0 \notin
\Omega$, the integral $\int_{1/2}^1 |v'(\rho)|^2 \rho^{N-1} \, d\rho$
is finite if and only if $\int_{1/2}^1 |v'(\rho)|^2 \, d\rho$ is
finite.
  Therefore,
\begin{equation} \label{eq:7}
R_{\a,p}(u) \geq {C}{\alpha^{1+\frac{2}{p}}} \inf_{\substack{ v
\in H^1_0 (\Omega) \\ v \neq 0}} \frac{\int_{1}^3 |v'(\rho)|^2 |4
\rho -3| \, d\rho}{\left(\int_{1}^3 |v(\rho)|^p \, d\rho
\right)^{2/p}}.
\end{equation}
where $C$ depends only on $N$. To complete the proof, we will show
that the right-hand side of \eqref{eq:7} is greater than zero.
This follows from some general Hardy--type inequality (see
\cite{opic}, Theorem 1.14), but we present here an elementary
proof for the sake of completeness. Indeed, given $v \in
H^1_{0,\mathrm{rad}} (\Omega)$, we can write for $\rho\in [1,2]$
\begin{eqnarray*}
|v(\rho)| &=& |v(\rho)-v(1)| \leq \int_{1}^\rho |v'(t)|
|2-t|^{1/2} \frac{dt}
{|2-t|^{1/2}} \\
&\leq& \left( \int_{1}^\rho |v'(t)|^2 |2-t| \, dt \right)^{1/2}
\left( \int_{1}^\rho \frac{dt}{|2-t|} \right)^{1/2} \\
&\leq& \left( \int_{1}^3 |v'(t)|^2 |2-t| \, dt \right)^{1/2}
\left( - \log |2-\rho| \right)^{1/2}.
\end{eqnarray*}
Hence
\begin{eqnarray*}
  \int_{1}^{2} |v(\rho)|^p \, d\rho &\leq&  \left( \int_{1}^3 |v'(\rho)|^2 |2-\rho| \, d\rho \right)^{p/2}
   \int_{1}^{2} \left( -\log (2-\rho) \right)^{p/2} \, d\rho \\
  &=&  \left( \int_{1}^3 |v'(\rho)|^2 |2-\rho| \, d\rho \right)^{p/2}
  \int_0^\infty t^{p/2} e^{-t} \, dt \\
  &=&  \Gamma \left( \frac{p+2}{2} \right)
  \left( \int_{1}^3 |v'(\rho)|^2 |2-\rho| \, d\rho \right)^{p/2},
\end{eqnarray*}
and in a similar way
\[
\int_{2}^{3} |v(\rho)|^p\, d\rho \leq  \Gamma \left( \frac{p+2}{2}
\right) \left( \int_{1}^3 |v'(\rho)|^2 |2-\rho| \, d\rho
\right)^{p/2}.
\]
Therefore
\[
{\int_{1}^{3} |v'(\rho)|^2 |2-\rho| \, d\rho} \geq \frac{1}{2^{
2/p}\Gamma \left( \frac{p+2}{2}
\right)^{2/p}}\left(\int_{1}^3|v(\rho)|^p\;d\rho\right)^{\frac
2p}.
\]
This implies that  the infimum in \eqref{eq:7} is strictly
positive and for any $M>2$ there exists a constant $C_1>0$ such
that $2^{-2/p}\geq C_1\, \Gamma \left( \frac{p+2}{2}
\right)^{2/p}$ for any $p\in (2,M]$, since the Gamma function is
positive, $C^\infty$ and  $\Gamma \left( \frac{p+2}{2} \right)\sim
(p/2)^{p/2}{\rm e}^{-p/2}\sqrt{\pi p}$ for $p\to +\infty$. We
finally collect \eqref{eq:8} and \eqref{eq:14} to get the desired
conclusion
\[
S_{\alpha,p}^{\mathrm{rad}} \geq C_1 \, \alpha^{1+\frac{2}{p}}.
\]
\end{proof}

Set now
\begin{equation}\label{eq3}
S_{\alpha,p}=\inf_{\substack{u \in H_0^1(\Omega) \\ u \neq
0}}R_{\alpha,p}(u).
\end{equation}
It is easily proved that for $p$ subcritical ~$S_{\alpha,p}$ is
attained by a function $u_{\alpha,p}$ that  satisfies (up to a
scaling) equation \eqref{eq1}.

In order to prove that { any} solution  $u_{\alpha,p}$ is not
radial (at least for $\a$ large) we need
 an estimate from above of the level
$S_{\a,p}$.
\begin{lemma}
Let $p\in (2,2^*)$. There exists $\bar \alpha$ such that for
$\a\geq\bar\alpha$
\begin{equation} \label{eq:S-asym}
S_{\a,p}\le C\a^{2-N+\frac{2N}{p}}.
\end{equation}
\end{lemma}
\begin{proof} The proof essentially follows the same techniques    of
\cite{ssw}.

Let $\psi$ be a positive smooth function with support in the unit
ball
  $B$. Let us consider
  $\psi_\a(x)=\psi(\a(x-x_\a))$, where $x_{\a}=(3-\frac
  1\a,0,\ldots,0)$. Since $\psi_\a$ has support in the ball
  $B(x_\a,\frac 1\a)$, by the change of variable $y=\a(x-x_\a)$ we
  obtain:
\begin{equation*}
\int_{\O} \Psi_\alpha(x) \psi^p_{\a}(x)\, dx =
\int_{B(x_\alpha,\frac 1\alpha)}||x|-2|^\alpha \psi_\a^p(x)\, dx
\geq\left(1-\frac 2\alpha\right)^\alpha \alpha^{-N}\int_{B}
\psi^p(y)\, dy
\end{equation*}
Moreover
\begin{equation*}
\int_{\O} |\nabla \psi_{\a}|^2\, dx = \a^2\int_{\O} |\nabla
\psi(\a(x-x_\a))|^2\, dx =\a^{2-N}\int_{B}|\nabla \psi|^2\, dx,
\end{equation*}
so that
\begin{equation*}
S_{\a,p}\le R_{\a,p}(\psi_\a)\le  C{\alpha^{2-N+\frac{2N}{p}}}
\end{equation*}
for all $\alpha$ sufficiently large. This proves the Lemma.
\end{proof}

By comparing \eqref{eq:S-asym} and \eqref{eq:11}, we deduce a
first symmetry--breaking result.

\begin{theorem}
  Let $p \in (2,2^*)$. For $\alpha$  sufficiently large, { any} ground
  state $u_{\alpha,p}$ is a non-radial function.
\end{theorem}
\begin{proof}
  From \eqref{eq:S-asym} and \eqref{eq:11} it follows that
  $S_{\alpha,p} < S_{\alpha,p}^{\mathrm{rad}}$ when $\alpha$ is large.
\end{proof}

\vspace{.2cm}
\section{Symmetry breaking as $p\to 2^*$}

In this section we consider  $\alpha$ fixed, $p$ close to $2^*$
and we establish the following
\begin{theorem}\label{3sol}
Let $\a>0$. For $p $ close to $2^*$ the quotient $R_{\alpha,p}$
has at least two non radial local minima.
\end{theorem}

We briefly explain how the proof proceeds. Of course, we already
know that any global minimizer of $R_{\alpha,p}$ produces a first
solution $u_{\alpha,p}$.  As the Theorem \ref{conc} states, this
solution concentrates at precisely one point of the boundary
$\partial\Omega$. Since this boundary has two connected
components, we will minimize $R_{\alpha,p}$ over the set $\Lambda$
of $H_0^1$ functions which ``concentrate'' at the opposite
component of the boundary. A careful estimate is proved in order
to show that minimizers fall inside the interior of $\Lambda$.

Consider now any minimizer~$u_{\alpha,p}$. As in \cite{caopeng} we have a
description of the profile of $u_{\alpha,p}$ as $p \to 2^*$.

\begin{theorem}\label{conc}
Let $p\in (2,2^*)$ and $\alpha >0$.  { Any} minimum $u_{\alpha,p}$
of $R_{\alpha, p}(u)$ in $H_0^1\setminus\{ 0\}$ satisfies (passing
to a subsequence) for some $x_0\in \partial \Omega$
\begin{itemize}
  \item[i)] $|\nabla u_{\alpha,p}|^2 \rightarrow \mu \delta_{x_0}$
   weakly in the sense of measure as $p\rightarrow 2^*$,
\item[ii)] $ |u_{\alpha,p}|^{2^*}\rightarrow \nu \delta_{x_0}$
weakly in the sense of measure as $p\rightarrow 2^*$,
\end{itemize}
where $\mu>0$ and $\nu >0$ are such that $\mu \geq
S_{0,2^*}\nu^{2/2^*}$ and $\delta_x$ is the Dirac mass at $x$.
\end{theorem}

\begin{proof}
This result can be proven by using, with suitable modifications,
the same arguments of \cite{caopeng}.
\end{proof}
 To get a second local minimizer, we will assume without loss
of generality that { any} $u_{\alpha,p}$ concentrates at some
point on the sphere $|x|=3$ (a similar argument holds if
$u_{\alpha,p}$ concentrates at some point $x$ with $|x|=1$). After
a rotation, we can even assume that { any} $u_{\alpha,p}$
concentrates at the point $(3,0,\dots,0)$.

Let
\[
\Omega^-=\left\{x\in\mathbb{R}^N \mid 1<|x|<2 \right\} ,\quad
  \Omega^+=\left\{x\in\mathbb{R}^N \mid 2<|x|<3 \right\}
\]
and
\[
\Sigma =\left\{u\in H_0^1 \setminus\{0\} \mid  \int_{\Omega
^+}|\nabla u|^2\ dx=
 \int_{\Omega ^-}|\nabla u|^2\ dx\right\}.
\]
Let us denote
\[
T_{\alpha,p}=\inf_{u \in \Sigma} R_{\alpha,p}(u).
\]
We have the following uniform estimate for $T_{\alpha,p}$.
\begin{proposition}\label{equipart}
  Let $\a>0$. There exists $\delta>0$ such that
\[
\liminf_{p\to 2^*}T_{\alpha,p}>S_{0,2^*}+\delta.
\]
\end{proposition}
\begin{proof}
We first prove that  $T_{\alpha,p} $ is achieved by a function
$v_{\a,p}\in \Sigma$. Consider a minimizing sequence $\{u_n\}$ for
$T_{\alpha,p} $. We can exploit the homogeneity of $R_{\alpha,p}$
and assume that $\int_{\Omega}|\nabla u_n|^2 \, dx=1 $. Up to a
subsequence, $u_n$ converges to
$v=v_{\a,p}$  weakly in $H^1_0(\Omega)$ and strongly in
$L^q(\Omega)$, for all $q \in (2,2^*)$. All we have to check is
that $v\in \Sigma$ (proving \textit{a posteriori} that the
convergence of $u_n$ to $v$ is strong). From the strong convergence in
$L^q(\Omega)$ we have that
\begin{equation} \label{sigma}
R_{\alpha,p}(v)\le\frac{1}{\left(\int_{\Omega}\Psi_\alpha(x)|v|^p
dx\right)^{2/p}}=T_{\alpha,p}
\end{equation}
and in particular $v\neq 0$. It is not restrictive to suppose that
$v \geq 0$ in $\Omega$. For the sake of contraddiction, assume that
$$
\int_{\Omega^+}|\nabla v|^2<\frac12.
$$
Fix a nonnegative smooth function $\varphi_1\in
C_0^{\infty}({\Omega^+})$,  $\varphi_1\neq 0$ and $\delta\geq 0$.
Setting $u=v+\delta\varphi_1$ from the positivity of $v$ and
$\varphi_1$, we have, for $\delta>0$,

\begin{equation} \label{sigma2}
\int_{\Omega}\Psi_\alpha(x)|v|^p
dx<\int_{\Omega}\Psi_\alpha(x)|u|^p dx.
\end{equation}

Now,
$$
\int_{\Omega^+}|\nabla u|^2 dx=\int_{\Omega^+}|\nabla
v|^2dx+2\delta\int_{\Omega^+}\nabla v\cdot\nabla\varphi_1 dx
+\delta^2\int_{\Omega^+}|\nabla \varphi_1|^2dx
$$
If we define $g_1:[0,+\infty)\to \mathbb R$ by
$$
g_1(\delta)=\int_{\Omega^+}|\nabla
v|^2dx+2\delta\int_{\Omega^+}\nabla v\cdot\nabla\varphi_1 dx
+\delta^2\int_{\Omega^+}|\nabla \varphi_1|^2dx
$$
we see that $g_1$ is continuous and $ g_1(0)<\frac
12,\quad\displaystyle \lim_{\delta\to+\infty}g_1(\delta)=+\infty.
$ Hence there exists $\delta_1>0$ with $g_1(\delta_1)=\frac12$. We
 can reason in an analogous way if $\int_{\Omega^-}|\nabla
v|^2dx< 1/2$ in order to find $\delta_2\geq 0$ and $ \varphi_2\geq
0$  such that $\int_{\Omega^-}|\nabla
(v+\delta_2\varphi_2)|^2dx=\frac 12$ .

 From  \eqref{sigma2},  this
shows that there exists
$w=v+\delta_1\varphi_1+\delta_2\varphi_2\in \Sigma$ such that
$$
R_{\alpha,p}(w)<T_{\alpha,p}
$$
which gives a contraddiction. Finally we must have that
$v_{\alpha,p}\in\Sigma,$  is a minimum point.

\vspace{.2cm}

 Moreover for any $\alpha>0,$ and $ 2<p<2^*$
we have
 $T_{\alpha,p}\geq S_{\alpha,p}$. We want to prove that  the inequality is
 strict  at least for $p\to 2^*$. Indeed assume on the contrary
 that
 \[
 \liminf_{p\to 2^*}T_{\alpha,p}=\liminf_{p\to
 2^*}R_{\alpha,p}(v_{\alpha,p})=S_{0,2^*}.
 \]
From the definition of $S_{0,2^*}$ and H$\ddot{\rm o}$lder
inequality we get, for a subsequence $p=p_k\to 2^*$
\begin{multline} \label{eq:comebrezis}
S_{0,2^*}\leq \frac{\int_{\Omega}|\nabla v_{\a,p}|^2\
dx}{\left(\int_{\Omega}| v_{\a,p}|^{2^*}\
dx\right)^{2/{2^*}}}\leq{
|\Omega|^{\frac{(2^*-p)2}{2^*p}}}\frac{\int_{\Omega}|\nabla
v_{\a,p}|^2\ dx}{(\int_{\Omega}| v_{\a,p}|^{p}\
dx)^{2/{p}}} \\
\leq{ |\Omega|^{\frac{(2^*-p)2}{2^*p}}}
 \frac{\int |\nabla
v_{\a,p}|^2\ dx}{(\int_{\Omega}\Psi_\alpha(x)| v_{\a,p}|^{p}\
dx)^{2/{p}}}=S_{0,2^*}+o(1)
\end{multline}
since the weight satisfies  $\Psi_\alpha(x)\leq 1$. In particular
\[
\frac{\int_{\Omega}|\nabla v_{\a,p}|^2\ dx}{\left(\int_{\Omega}|
v_{\a,p}|^{2^*}\ dx\right)^{2/{2^*}}}\to S_{0,2^*},
\]
and $v_{\a,p}$ is a minimizing sequence of $S_{0,2^*}$.

In the same way as Cao and Peng did in {\cite{caopeng}}, Theorem
1.1,  we can prove that $v_{\alpha,p}$ concentrates at precisely
one point one of the boundary $\partial \Omega$. This contradicts
the fact that $\int_{\Omega ^+}|\nabla v_{\alpha,p}|^2\ dx=
 \int_{\Omega ^-}|\nabla v_{\alpha,p}|^2\
dx$.
\end{proof}

Consider now the points
\[
x_{0,\varepsilon} = x_0= \left(3-\frac{1}{|\log
\varepsilon|},0,\ldots,0\right), \quad x_{1,\varepsilon} = x_1=
\left( 1+\frac{1}{|\log \varepsilon|},0,\ldots,0 \right)
\]
and
\[
U(x) = \frac{1}{\left(1 +|x| \right)^{(N-2)/2}}.
\]
We recall that $S_{0,2^*}$ is not achieved on any proper subset of
$\R^N$, and that it is independent of $\Omega$. However, it is known
that $S_{0,2^*}(\R^N)$ is achieved, and all the minimizers can be
written in the form
\[
\mathcal{U}_{\theta,y}(x)=\frac{1}{\left(\theta^2+|x-y|^2
  \right)^\frac{N-2}{2}}, \quad \hbox{$\theta>0$, $y \in \R^N$}.
\]

We set
\[
U^i_\varepsilon(x)=\varepsilon^{-\frac{N-2}{2}} \
U\left(\frac{x-x_i}{\sqrt\varepsilon}\right)=\frac{1}{(\varepsilon+|x-x_i|^2)^{\frac{N-2}{2}}},
\]
and denote by $\varphi_i$ ($i=0,1$) two \textit{cut-off} functions
such that $0\le\varphi_i\le 1$, $|\nabla\varphi_i|\le C|\log
\varepsilon|$ for some constant $C>0$, and
\[
\varphi_i (x) =
\begin{cases}
  1, &\text{if $|x-x_i| < \frac{1}{2|\log \varepsilon|}$}  \\[10pt]
  0, &\text{if $|x-x_i| \geq \frac{1}{|\log \varepsilon|}$}.
\end{cases}
\]
The following Lemma shows that the truncated functions
\begin{equation} \label{eq:instanton}
  u^i_\varepsilon(x)=\varphi_i(x)U^i_\varepsilon(x), \quad i=0,1,
\end{equation}
are almost minimizers for $S_{0,2^*}$. We omit the proof of this fact,
since it is an easy modification of the argument of Cao and Peng in
\cite{caopeng}.

\begin {lemma}\label{zero}
  Let $\a >0$. There results
  \[
  \lim_{p\to
    2^*}R_{\alpha,p}(u^i_\varepsilon)=S_{0,2^*}+K(\varepsilon)
  \]
  with $\displaystyle\lim_{\varepsilon \to 0}K(\varepsilon)=0.$
\end{lemma}

\begin{remark}
A direct consequence of Lemma \ref{zero} is that
    $S_{0,2^*}=S_{\alpha, 2^*}$. Indeed $S_{0,2^*}\leq S_{\alpha,2^*}$
    since $\Psi_\alpha (|x|)\leq 1$ . On the other hand by Lemma
    \ref{zero}
\[
 R_{\alpha,2^*}(u^i_\varepsilon)= \lim_{p\to
 2^*}R_{\alpha,p}(u^i_\varepsilon)=S_{0,2^*}+K(\varepsilon).
\]
Therefore $S_{0,2^*}+K(\varepsilon)\geq S_{\alpha,2^*}$ for every
$\varepsilon >0$. Letting $\varepsilon \to 0$ we conclude
$S_{0,2^*}\geq S_{\alpha,2^*}$.
\end{remark}

\bigskip

\noindent We are now ready to prove the Theorem \ref{3sol}.
\begin{proof}[Proof of Theorem \ref{3sol}]
 Let $u_{\alpha,p}$ be a ground state solution. Let
us suppose that it concentrates on the outer boundary.  Consider
the open subset
\[
\Lambda=\left\{u\in H^1_0(\Omega):\;\int_{\Omega^-}|\nabla u|^2\
dx> \int_{\Omega^+}|\nabla u|^2\ dx\right\}.
\]
The infimum of $R_{\alpha,p}$ on $\overline\Lambda$ is achieved.
However it cannot be achieved on the boundary $\partial
\Lambda=\Sigma$. Indeed, by Proposition \ref{equipart},
$$\displaystyle \inf_\Sigma
R_{\alpha,p}>S_{0,2^*}+\delta \;\;{\rm as }\;\;p\to 2^*
$$
 and $$\displaystyle\inf_{\Lambda}R_{\alpha,p}(u)\le
R_{\alpha,p}(u^1_{\varepsilon})\to
S_{0,2^*}+K_1(\varepsilon)\;\;{\rm as }\;\;p\to 2^*$$
 since
$u^1_\varepsilon\in \Lambda $ for $\varepsilon $ small enough.
Then the infimum is achieved in a interior point of $\Lambda$ and
is therefore a critical point of $R_{\alpha,p}$.

\end{proof}

\begin{remark}
  Theorem \ref{conc} shows that { any} ground state solution
 ``concentrates" in a single point at the boundary as $p\to 2 ^*$ and
  consequently this solution is not radial. This symmetry breaking can
  be also proved by using a continuation argument as in
  \cite{BN}. Indeed, \eqref{eq:comebrezis} shows that $\lim_{p \to
    2^*} S_{\alpha,p} = S_{0,2^*}$, and since $S_{0,2^*}<S^{\rm
    rad}_{0,2^*}$ we conclude as in \cite{BN} that ground states of
  $S_{\alpha,p}$ cannot be radially symmetric as $p \to 2^*$.
\end{remark}

\section{existence of a third non-radial solution}

In the previous section we proved the existence of two solutions
of \eqref{eq1} which are  local minima of the Rayleigh quotient
for $p$ near $2^*$. One would expect another critical point of
$R_{\alpha,p}$ located in some sense between these minimum points.
This is precisely the idea we are going to pursue further in the
current section.

For ${\varepsilon}$ small enough let $u^i_{\varepsilon}=\varphi_i
U^i_{\varepsilon}$, $i\in \{0,1\}$, be defined as in
\eqref{eq:instanton}.  We will prove that $R_{\alpha,p}$ has the
Mountain Pass geometry.

\bigskip

Let us introduce the \textit{mountain--pass level}
\[
\beta=\beta(\alpha,p)=\inf_{\gamma\in\Gamma}\max_{t\in
[0,1]}R_{\alpha,p}(\gamma(t)),
\]
where $\Gamma =\{\gamma\in C([0,1],H^1_0(\Omega)) \mid
\gamma(0)=u^0_{\varepsilon},\; \gamma(1)=u^1_{\varepsilon}\}$ is
the set of continuous paths joining $u_\ge^0$ with $u_\ge^1$.
We claim that $\beta$ is a critical value for $R_{\alpha,p}$.

We begin to prove that $\beta$ is larger, uniformly with respect
to~$\ge$, than the values of the functional $R_{\alpha,p}$ at the
points $u_\ge^0$ and $u_\ge^1$.

\begin{lemma}\label{lemma:cont} Let
  $M_{\varepsilon}=\max\{R_{\alpha,p}(u^0_{\varepsilon}),R_{\alpha,p}(u^1_{\varepsilon})\}$.
  There exists $\sigma>0$ such that $\beta\geq
  M_{\varepsilon}+\sigma$ uniformly with respect
to~$\ge$.
\end{lemma}

\begin{proof}
  We prove that there exists $\sigma$ such that for all $\gamma \in
  \Gamma$
\[
\max R_{\alpha,p}(\gamma(t))\geq M_{\varepsilon}+\sigma.
\]
A simple continuity argument shows that for every $\gamma \in
\Gamma$ there exists $t_{\gamma}$ such that $\gamma(t_{\gamma})\in
\Sigma$, where
\[
\Sigma =\left\{u\in H_0^1 \setminus\{0\} \mid  \int_{\Omega
^+}|\nabla u|^2\ dx=
 \int_{\Omega ^-}|\nabla u|^2\
dx\right\}.
\]
Indeed the map $t \in [0,1] \mapsto \int_{\Omega^{+}} |\nabla
\gamma (t)|^2 \, dx - \int_{\Omega^{-}} |\nabla \gamma (t)|^2 \,
dx$ is continuous and it  takes a negative value at $t=0$ and a
positive value at $t=1$. Now Proposition \ref{equipart} implies,
for $p$ near $2^*$ the existence of $\delta>0$ with
\[
\max_{t\in [0,1]} R_{\alpha,p}(\gamma(t))\geq
R_{\alpha,p}(\gamma(t_{\gamma}))\geq \inf_{u\in \Sigma
}R_{\alpha,p}(u)\geq S_{0,2^*}+\delta.
\]
On the other hand, for $\varepsilon$ sufficiently small,
\[
M_{\varepsilon}<S_{0,2^*}+\frac{\delta}{2}.
\]
This concludes the proof.
\end{proof}

The previous estimate allows us to show that $\beta$ is a critical
level for $R_{\alpha,p}$. Therefore a further nonradial solution
to \eqref{eq1} arises.

\begin{proposition}\label{prop:MP}
There exist  $\bar{\alpha} > 0$ and $2<\bar p<2^*$ such that   for
all $\a\geq\bar{\a}$ and $\bar p\le p< 2^*$ it results that
$\beta$ is a critical value for $R_{\a,p}$ and it is attained by a
non-radial function.
\end{proposition}
\begin{proof}
  From the previous result we can apply a deformation argument (see
  \cite{AM,st}) to prove that $\beta$ is a critical level and it is
  attained (since the PS condition is satisfied) by a
  function $w$. From the asymptotic estimate \eqref{eq:11} for the
  radial level $S_{\alpha,p}^{\mathrm{rad}}$, one has that there exists a
  constant $C$ independent from $p$ such that
  \[
  S_{\alpha,p}^{\mathrm{rad}}\geq C \alpha^{1+2/p}.
  \]
  In particular
  \[
  S_{\alpha,p}^{\mathrm{rad}}\to +\infty \quad\text{as $\alpha \to
    +\infty$}.
  \]
  This allows us to choose $\alpha_0$ such that
  $S_{\alpha,p}^{\mathrm{rad}}\geq 3 S_{0,2^*}$ for all
  $\alpha\geq\alpha_0$.

  Define $\zeta\in\Gamma$ by $\zeta(t)=t
  u^1_{\varepsilon}+(1-t)u^0_{\varepsilon}$ for all~$t \in [0,1]$, and
  let $\tau \in [0,1]$ be such that $R_{\alpha,p}(\zeta(
  \tau))=\max_{t\in [0,1]} R_{\alpha,p}(\zeta(t))$.

\vspace{.1cm}

Since $u^1_{\varepsilon}$ and $u^0_{\varepsilon}$ have disjoint
supports one has, for $\varepsilon $ sufficiently small,
\begin{eqnarray*}
  R_{\alpha,p}(w) &=& \beta\le R_{\alpha,p}(\zeta(\tau))=
  \frac{\int_{\Omega}|\nabla (\tau
    u^1_{\varepsilon}+(1-\tau)u^0_{\varepsilon})|^2\ dx}
  {\left(\int_{\Omega}\Psi_{\alpha}|\tau
      u^1_{\varepsilon}+(1-\tau)u^0_{\varepsilon}|^p \ dx\right)^{2/p}}
  \\[10pt]
  &=&
  \frac{\int_{\Omega}\tau^2|\nabla
    u^1_{\varepsilon}|^2\ dx+ \int_{\Omega}(1-\tau)^2|\nabla
    u^0_{\varepsilon}|^2\ dx}{{\left(\tau^p \int_{\Omega}\Psi_{\alpha}|
        u^1_{\varepsilon}|^p\
        dx+(1-\tau)^p\int_{\Omega}\Psi_{\alpha}|u^0_{\varepsilon}|^p \ dx\right)^{2/p}}}
  \\[10pt]
  &\leq& \frac{\tau^2\int_{\Omega}|\nabla u^1_{\varepsilon}|^2\
    dx}{{\left(\tau^p \int_{\Omega}\Psi_{\alpha}|
        u^1_{\varepsilon}|^p\ dx\right)^{2/p}}}+\frac{(1-\tau)^2\int_{\Omega}|\nabla
    u^0_{\varepsilon}|^2\ dx}{{\left((1-\tau)^p
        \int_{\Omega}\Psi_{\alpha}|
        u^0_{\varepsilon}|^p\ dx\right)^{2/p}}} \\[10pt]
  &=&R_{\alpha,p}(u^0_{\varepsilon})+R_{\alpha,p}(u^1_{\varepsilon})\leq
  2 M_{\varepsilon}<3 S_{0,2^*}\le S_{\alpha,p}^{\mathrm{rad}}.
\end{eqnarray*}
This concludes the proof.
\end{proof}
\vspace{.2cm}

\section{Behaviour of the ground-state solutions for $\alpha$
large}

This section is devoted to the analysis of  a ground state
solution as $\alpha \to +\infty$. Even in this case this solution
tends to ``concentrate" at the boundary~$\partial\Omega$. However,
this concentration is much weaker than the concentration as $p \to
2^*$.

\bigskip

We use the notation~$C(r_1,r_2)=\{x\in \mathbb R^N \mid
r_1<|x|<r_2\}$. Let $\delta$ be sufficiently small (say
$\delta<\frac12$) and $\phi$ be a smooth cut-off function such
that $0\le\phi\le 1$ with
\begin{equation} \label{eq:phi}
\phi(x) =
\begin{cases}
1, &x\in C(1,1+\delta) \cup C(3-\delta,3) \\
0, &x\in C(2-\delta,2+\delta)
\end{cases}
\end{equation}
From now on, since $p\in(2,2^*)$ is fixed we denote a ground state
solution of problem (\ref{eq1}) $u_{\a,p}$ with $u_{\a}$.
\begin{proposition}\label{prop:1}
  Let $u_{\alpha}$ be such that $R_{\alpha,p}(u_\a)=S_{\alpha,p}$. If $\phi$ is
  as in \eqref{eq:phi}, then
\begin{equation}
R_{\alpha, p}(\phi u_\a) = S_{\alpha,p} + o(S_{\alpha,p})
\quad\hbox{as $\alpha \to +\infty$}.
\end{equation}
\end{proposition}
\begin{proof}
It is not restrictive, by the homogeneity of $R_{\alpha,p}$, to
assume $\int_{\Omega}|\nabla u_\a|^2\ dx=1$. We split the proof
into two steps.

\noindent\textbf{Step 1.} We claim that
\begin{equation} \label{eq:1}
\int_{\Omega}\psa (u_\a \phi)^p\ dx= \int_{\Omega}\psa u_\a^p\ dx\
+ o \left( \int_{\Omega}\psa u_\a^p\ dx \right)
\end{equation}
Indeed, suppose on the contrary that
\[
\limsup_{\alpha\to\infty}\frac{\int_{\Omega}\psa u_\a^p(1-\phi^p)\
dx}{\int_{\Omega}\psa u_\a^p\ dx}=\beta>0
\]
This implies that, up to some subsequence,
\[
\frac{\int_{\Omega}\psa u_\a^p(1-\phi^p)\ dx}{\int_{\Omega}\psa
u_\a^p\ dx}>\beta/2>0
\]
Since $1-\phi^p\equiv 0$ on $C(1,1 +\delta) \cup C(3-\delta,3)$ we
have
\begin{eqnarray*}
\int_{\Omega}\psa u_\a^p(1-\phi^p)\ dx &=&
\int_{C{(1+\delta,3-\delta)}}\psa u_\a^p(1-\phi^p)\
dx \\
&\leq& (1-\delta)^{\alpha}\int_{\Omega}u_{\a}^p(1-\phi^p)\
dx\le(1-\delta)^{\alpha}\int_{\Omega}u_{\a}^p\ dx.
\end{eqnarray*}
Therefore
\[
\int_{\Omega}u_{\a}^p\ dx\geq
(1-\delta)^{-\alpha}\int_{\Omega}\psa u_\a^p(1-\phi^p)\ dx
\]
Now
\[
\frac{\int_{\Omega}u_{\a}^p\ dx}{\int_{\Omega}\psa u_{\a}^p\
dx}\geq (1-\delta)^{-\alpha}\, \frac{\int_{\Omega}\psa
u_\a^p(1-\phi^p)\ dx}{{\int_{\Omega}\psa u_{\a}^p}\ dx}\geq
(1-\delta)^{-\alpha}\, {\frac{\beta}{2}}.
\]
Since $S_{\alpha,p}^{p/2}=\left( \int_{\Omega}\psa u_\a^p\, dx
\right)^{-1}$ the last inequality can be written as
\[
S_{\alpha,p}^{p/2}\geq \frac{\beta}{2} \,
\frac{(1-\delta)^{-\alpha}}{{\int_{\Omega}u_{\a}^p}\ dx} \geq
{{\frac{\beta}{2}}} (1-\delta)^{-\alpha}S^{p/2}_{0 ,p},
\]
where
\[
S_{0,p}=\inf_{u\neq 0}\frac{\int_{\Omega}|\nabla u|^2\,
dx}{(\int_{\Omega} u^p\, dx)^{2/p}}
\]
On the other hand  from \eqref{eq:S-asym} one has the estimate
\[
S_{\alpha,p}^{p/2}\le C \alpha ^{p-\frac{N}{2}p+ N},
\]
which gives a contradiction for $\alpha$ large. This proves
\eqref{eq:1}.

\noindent \textbf{Step 2.} Now we prove that
\begin{equation} \label{eq:2}
\int_{\Omega} |\nabla {u_\a\phi}|^2\ dx= \int_{\Omega} |\nabla
{u_\a}|^2\ dx\; +\ o( 1) =1+o(1).
\end{equation}
It is not difficult to prove that $u_\a$ satisfies the problem
\begin{equation}\label{eq5}
\begin{cases}
 -\Delta u_\a=S_{\a,p}^{p/2}\Psi_{\alpha}u_\a^{p-1} &\text{in $\Omega$},\\
 u_\a>0 &\text{in $\Omega$},\\
u_\a =0 &\text{on $\partial \Omega$},
\end{cases}
\end{equation}
Since $\| \nabla u_\a \|_2=1$, up to subsequences, we have that,
as $\alpha \to \infty$,
\[
u_\alpha \to u \quad\hbox{weakly in $H_0^1(\Omega)$, strongly in
  $L^q(\Omega)$, and a.e.}
\]
We now prove that $u=0$. Indeed, multiplying equation \eqref{eq5}
by a smooth function $\psi$ with $\supp \psi  \subset \subset \O$
and integrate, we obtain
\[
\int_{\Omega} \nabla u_\a \nabla\psi\, dx=\int_{\Omega}
S_{\a,p}^{p/2} \Psi_{\alpha}u_\a^{p-1}\psi\, dx \to 0 , \quad
\alpha\to+\infty
\]
since, by \eqref{eq:S-asym},  $S_{\alpha,p}^{p/2}\psa\to 0$
uniformly on $\supp \psi$ and $u_\a$ is uniformly bounded in $L^q$
for $1\le q< 2^*$. Hence $\int_{\Omega} \nabla u \cdot \nabla
\varphi \, dx = 0$ for all $\varphi \in C_0^\infty(\Omega)$. Since
$u \in H_0^1(\Omega)$, this implies that $u=0$.

Now we estimate the difference
\begin{multline} \label{eq:20}
\left| \int_{\Omega} |\nabla {u_\a}|^2\;dx-\int_{\Omega} |\nabla
{(u_\a\phi)}|^2\;dx\right|\le
\\
\le \int_{\Omega} |\nabla {u_\a}|^2(1-\phi^2)\;dx+\int_{\Omega}
|\nabla \phi|^2u_\a^2\;dx+2\left|\int_{\Omega} \nabla {u_\a}\nabla
\phi {u_\a}\phi\;dx\right|
\end{multline}
The last terms tend to zero thanks to the strong convergence in
$L^q$ for all $q\in[1,2^*)$.  In order to estimate the term
$\int_{\Omega} |\nabla {u_\a}|^2(1-\phi^2)\ dx$, we multiply
\eqref{eq5} by $u_\a(1-\phi^2)=u_\a\eta$ and integrate.  Since
$\supp\eta=\supp (1-\phi^2)\subset\subset\Omega$ we have
\[
\int_{\Omega} \nabla {u_\a}\nabla(\eta u_\a )\,
dx=\int_{\Omega}S_{\a,p}^{p/2}\Psi_{\alpha}u_\a^{p}\eta\, dx ,
\]
so that
\begin{eqnarray*}
\left|\int_{\Omega}| \nabla {u_\a}|^2\eta\, dx\right| &\leq&
\left|\int_{\Omega}u_\a\nabla\eta\nabla u_\a\, dx\right|+
\left|\int_{\Omega}S_{\a,p}^{p/2}\Psi_{\alpha}u_\a^{p}\eta\, dx
\right|
\\
&\leq& ||\nabla\eta||_{\infty}\int_{\supp \eta}|\nabla u_\a u_\a\,
| \, dx+\left|\int_{\supp
\eta}S_{\a,p}^{p/2}\Psi_{\alpha}u_\a^{p}\eta
 \, dx \right|\to 0.
\end{eqnarray*}
\end{proof}

In proposition \ref{prop:1} we proved that the infimum of the
Rayleigh quotient $R_{\alpha,p}$ is essentially attained by the
function $\phi u_\alpha$. Thanks to the definition of $\phi$, we
can decompose $\phi u_\alpha=u_{\alpha,
  1}+u_{\alpha, 2}$, where $u_{\alpha, 1}$ vanishes in $C(2-\delta,
3)$ and $u_{\alpha, 2}$ vanishes in $C( 1,2 +\delta )$. The
following proposition is the main  step in order to prove that the
function $\phi u_\alpha$ concentrates at  the boundary.

\begin{proposition} \label{prop:3}
Let $\phi u_\alpha=u_{\alpha, 1}+u_{\alpha, 2}$, where $\operatorname
 {supp} u_{\alpha, 1} \subset C(1,2 -\delta)$ and
 $\operatorname{supp}u_{\alpha, 2} \subset C(2 +\delta,3)$, and
 $\lambda_\alpha = {\int_\Omega \Psi_\alpha u_{\alpha,1}^p \, dx}
 \big/ {\int_\Omega \Psi_\alpha u_{\alpha,2}^p \, dx}$. If
 $\lim_{n\to \infty}\lambda_{\alpha_n}=L$ for a sequence
 $\alpha_n\to \infty$ then either $L=0$ or $L=+\infty$.
\end{proposition}

\begin{remark}
For the quantity $\lambda_\alpha = {\int_\Omega \Psi_\alpha
u_{\alpha,1}^p \, dx} \big/ {\int_\Omega \Psi_\alpha
u_{\alpha,2}^p \, dx}$, we cannot exclude the case
$\limsup_{\alpha \to +\infty} \lambda_\alpha=+\infty$ and
$\liminf_{\alpha \to +\infty} \lambda_\alpha=0$. If  a
\emph{uniqueness} result for the minimizer $u_\alpha$ were known,
then it would be easy to conclude that $\alpha \mapsto u_\alpha$
is continuous. Therefore $\lambda_\alpha$ would be continuous,
too, and we could replace both the lower and the upper limit by
a unique limit. In general, one does not expect such a uniqueness
property for any $p$ and any $\alpha$.
 However, when $p \approx 2^*$ we suspect that the uniqueness argument of \cite{PiSe05} may be applied to our
 setting.
\end{remark}


\begin{proof}
By the definition of $u_{\alpha, 1}$ and $u_{\alpha, 2}$ we have
\begin{equation}
R_{\a,p}(\phi u_\alpha)=\frac{\int_{\Omega}|\nabla
u_{\alpha,1}|^2\, dx+\int_{\Omega}|\nabla u_{\alpha,2}|^2\,
dx}{\left(\int_{\Omega}\psa u_{\alpha,1}^p dx+\int_{\Omega}\psa
u_{\alpha,2}^pdx\right)^{\frac 2p}}.
\end{equation}
 Since
$u_\alpha$ is a positive solution, we can say that $\lambda_\alpha > 0$. We obtain the
following identity:
\begin{eqnarray}\label{eq7}
R_{\a,p}(\phi u_\alpha) &=&\frac{\int_{\Omega}|\nabla
u_{\alpha,1}|^2\, dx+\int_{\Omega}|\nabla u_{\alpha,2}|^2\,
dx}{\left(\lambda_\alpha\int_{\Omega}\psa
u_{\alpha,2}^p\, dx+\int_{\Omega}\psa u_{\alpha,2}^p\, dx\right)^{2/p}}\nonumber \\
&=& \frac{\int_{\Omega}|\nabla u_{\alpha,1}|^2\,
dx}{\left(\lambda_\alpha+1\right)^{2/p} \left(\int_{\Omega}\psa
u_{\alpha,2}^p \, dx\right)^{2/p}}+\frac{\int_{\Omega}|\nabla
u_{\alpha,2}|^2 \, dx}{\left(\lambda_\alpha+1\right)^{2/p}
\left(\int_{\Omega}\psa u_{\alpha,2}^p \, dx\right)^{2/p}} \nonumber \\
&=& \frac{\lambda_\alpha^{\frac 2p}\int_{\Omega}|\nabla
u_{\alpha,1}|^2\, dx}{\left(\lambda_\alpha+1\right)^\frac
2p\left(\int_{\Omega}\psa u_{\alpha,1}^p \,
dx\right)^{2/p}}+\frac{\int_{\Omega}|\nabla u_{\alpha,2}|^2\,
 dx}{\left(\lambda_\alpha+1\right)^\frac 2p\left(\int_{\Omega}\psa
u_{\alpha,2}^p \, dx\right)^{2/p}}.
\end{eqnarray}
By the definition of $S_{\alpha,p}$ each quotient
$R_{\a,p}(u_{\alpha,1})$ and $R_{\a,p}(u_{\alpha,2})$ in the last
term is greater than or equal to $S_{\alpha,p}$. Therefore  by
Proposition \ref{prop:1} and equation \eqref{eq7} one obtains
\begin{equation}\label{eq6}
S_{\alpha,p}+o(S_{\alpha,p})\geq \frac {1+\lambda_\alpha^\frac
2p}{(\lambda_\alpha+1)^\frac 2p}S_{\alpha,p}.
\end{equation}
We notice that the function $f(x)=\frac {1+x^{2/p}}{(x+1)^{2/p}}$
is strictly greater than 1 for every $x>0$, $f(0)=1$ and $f(x) \to
1$ as $x\to
  +\infty$. Moreover it is increasing in $[0,1]$ and
decreasing in $[1,+\infty)$ and $\max_{x >0} f(x) =
f(1)=2^{1-2/p}$. Let  $L\in \Lambda$ and  $\{\alpha_n\}$  a
sequence such that $\lambda_{\alpha_n} \to L$ as $n \to +\infty$.
Passing to the limit
 in \eqref{eq6}, we obtain that $1\geq\frac
{1+L^{2/p}}{(L+1)^{2/p}}$ and so either $L=+\infty$, or $L=0$.
\end{proof}

\begin{corollary} \label{coroll5}
With the notation of Proposition \ref{prop:3}, for any sequence $\{\alpha_n\}$
such that  $\lambda_{\alpha_n}\rightarrow 0$ one has
\begin{equation}
\lim_{n \to +\infty} \frac{\int_\Omega
|\nabla u_{\alpha_n,1} |^2\, dx}{\int_\Omega |\nabla
u_{\alpha_n,2}|^2 \, dx} =0.
\end{equation}
\end{corollary}
\begin{proof}
Let
\[
\xi_{\alpha}=\frac{\int_\Omega |\nabla u_{\alpha,1} |^2\,
dx}{\int_\Omega |\nabla u_{\alpha,2}|^2 \, dx}
\]
and  suppose that
$\limsup_{n\to\infty}\xi_{\a_n}> 0$. Up to subsequences,
$\xi_{\a_n}>\xi>0$ for some $\xi$. Therefore we have
\begin{eqnarray*}
S_{\alpha_n,p}+o(S_{\alpha_n,p}) &=& \frac{\int_{\Omega}|\nabla
u_{\alpha_n,1}|^2\ dx+\int_{\Omega}|\nabla u_{\alpha_n,2}|^2\
dx}{\left(\int_{\Omega} {\Psi_{\alpha_n}}u_{\alpha_n,1}^p
dx+\int_{\Omega}{\Psi_{\alpha_n}} u_{\alpha_n,2}^pdx\right)^{\frac
2p}}  \\
&=& \frac{(1+\xi_{\a_n})\int_{\Omega}|\nabla u_{\alpha_n,2}|^2\
dx}{\left(\int_{\Omega}{\Psi_{\alpha_n}}
u_{\alpha_n,2}^p dx\right)^{\frac 2p}(1+\lambda_{\a_n})} \\
&\geq& R_{\a_n,p}(u_{\a_n,2})\frac{1+\xi}{1+o(1)}\geq
(1+\xi)S_{\alpha_n,p}+o(S_{\alpha_n,p}),
\end{eqnarray*}
which is a contradiction. Hence
\[
\xi_{\alpha_n}=\frac{\int_\Omega |\nabla u_{\alpha_n,1} |^2\,
dx}{\int_\Omega |\nabla u_{\alpha_n,2}|^2 \, dx}\to 0.
\]
\end{proof}

\begin{remark}
An immediate consequence of the previous results is that in
particular for any $\alpha_n$ such that $\lambda_{\alpha_n}\to 0$
\begin{equation} \label{eq:strong}
\lim_{n \to +\infty} \int_\Omega |\nabla u_{\alpha_n,1}|^2 \, dx
=0.
\end{equation}
\end{remark}
\vspace{.2cm}

\section{Acknowledgements}
We would like to thank E.~Serra for suggesting the problem and for
his constant support.


\end{document}